\documentclass[12pt]{article}
\usepackage{graphics,amsmath,amssymb}
\usepackage{latexsym}
\usepackage{epsfig}
\usepackage{hyperref}
\usepackage{fancybox}

\def\epsilon{\varepsilon}

\def\phi{\varphi}

\newtheorem{theorem}{Theorem}[section]
\newtheorem{lemma}[theorem]{Lemma}

\newtheorem{definition}[theorem]{Definition}

\newtheorem{remark}[theorem]{Remark}

\def\N{{\mathbb N}}
\def\C{{\mathbb C}}
\def\R{{\mathbb R}}

\newenvironment{Proof}{\removelastskip\par\medskip
\noindent{\em Proof.} \rm}{\penalty-20\null$\square$\par\medbreak}

\title{\bf Fractional diffusion-wave equation: 
\\
hidden regularity for weak solutions}

\author{Paola Loreti
\thanks{Dipartimento di Scienze di Base e Applicate per l'Ingegneria,
Sapienza Universit\`a di Roma,
Via Antonio Scarpa 16, 00161 Roma (Italy); e-mail: 
$<$paola.loreti@sbai.uniroma1.it$>$ }
\and Daniela Sforza
\thanks{Dipartimento di Scienze di Base e Applicate per l'Ingegneria, 
Sapienza Universit\`a di Roma,
Via Antonio Scarpa 16, 00161 Roma (Italy); e-mail: 
$<$daniela.sforza@sbai.uniroma1.it$>$ }
\thanks{
This paper is published (in revised form) in Fract. Calc. Appl. Anal. Vol. 24, No 4 (2021), pp. 1015-1034,  DOI: 10.1515/fca-2021, and is available online at https://www.degruyter.com/journal/key/FCA/html
}}

\begin{document}

\maketitle

%
%
%

\bigskip
\noindent
\section{Introduction}
The equation 
\begin{equation}
\partial_t^{\alpha}u(t,x)=\triangle u(t,x)
\end{equation}
is obtained from the diffusion equation by replacing the first order time-derivative by the Caputo fractional derivative of order $\alpha$, where $1<\alpha<2$, that is
\begin{equation*}
\partial_t^{\alpha}f(t)=
\frac1{\Gamma(2-\alpha)}\int_0^t (t-\tau)^{1-\alpha}\frac{d^2f}{d\tau^2}(\tau)\ d\tau\,.
\end{equation*}
We will prove the following regularity results.
\begin{theorem}
If $u_0\in H^1_0(\Omega)$ and $u_1\in L^2(\Omega)$, then the unique weak solution $u$ of problem
\begin{equation}\label{eq:cauchy10}
\begin{cases}
\displaystyle
\partial_t^{\alpha}u(t,x) =\triangle u (t,x)\,,
\quad t\ge0,\,\, x\in \Omega,
\\
u(t,x)=0\qquad   t\ge0, \,\, x\in\partial\Omega,
\\
u(0,x)=u_{0}(x),\quad
u_t(0,x)=u_{1}(x),\qquad  x\in \Omega,
\end{cases}
\end{equation}
belongs to $C([0,T];H^1_0(\Omega))\cap C^1([0,T];D(A^{-\theta}))$, $\theta\in\big(\frac{2-\alpha}{2\alpha},\frac12\big]$, 
and
\begin{equation}\label{eq:in-data0}
\begin{split}
&\lim_{t\to0}\|u(t,\cdot)-u_0\|_{H^1_0(\Omega)}=
\lim_{t\to0}\|\partial_t u(t,\cdot)-u_1\|_{D(A^{-\theta})}=0\,,
\\
&\|u\|_{C([0,T];H^1_0(\Omega))}+\|\partial_t u\|_{C([0,T];D(A^{-\theta}))}
\le
C(\| u_0\|_{H^1_0(\Omega)}+\| u_1\|_{L^2(\Omega)}).
\end{split}
\end{equation}
In addition, for any $\theta\in\big(0,\frac1{2\alpha}\big)$ there exists a constant $C>0$ such that
\begin{equation}\label{eq:nabla-l20}
\|\nabla u\|_{L^2(0,T;D(A^{\theta}))}\le C\big(\|u_0\|_{H^1_0(\Omega)}+\|u_1\|_{L^2(\Omega)}\big),
\end{equation}
and for any $\theta\in\big(\frac{\alpha-1}{2\alpha},\frac12\big)$ there exists a constant $C>0$ such that
\begin{equation}\label{eq:partial-l20}
\|\partial_t^{\alpha}u\|_{L^2(0,T;D(A^{-\theta}))}\le C\big(\|u_0\|_{H^1_0(\Omega)}+\|u_1\|_{L^2(\Omega)}\big)\,.
\end{equation}
\end{theorem}

\begin{theorem}
Let  $u_0\in
H^1_0(\Omega)$ and $u_1\in L^{2}(\Omega)$. If $u$ is the weak solution of \eqref{eq:cauchy10}
then, for any $T>0$ there is a constant $c_0=c_0(T)$ such that, denoting by  $\partial_\nu u$ the normal derivative of $u$, we have
\begin{equation}\label{eq:hidden-alpha0}
\int_0^T\int_{\partial\Omega} \big|\partial_\nu u\big|^2d\sigma dt
\le c_0(\|\nabla u_0\|^2_{L^2(\Omega)}+\|u_1\|^2_{L^2(\Omega)})
\,.
\end{equation}

\end{theorem}

For previous results related to this problem see \cite {LasTri,Lio,MM,LoretiSforza2,LoretiSforza3} and references therein.


\section{Preliminaries}
Let $\Omega\subset\R^N$,  $N\ge1$, be a bounded open set with $C^2$ boundary. We consider 
$
L^2(\Omega)
$
endowed with the usual inner product and norm 
\begin{equation*}
\langle u,v\rangle=\int_{\Omega}u(x)v(x)\ dx,
\qquad
\|u\|_{L^2(\Omega)}=\left(\int_{\Omega}|u(x)|^{2}\ dx\right)^{1/2}\qquad
u,v\in L^2(\Omega)\,.
\end{equation*}
\begin{definition}
For any $f\in L^1(0,T)$ $(T>0)$ 
we define the Riemann--Liouville fractional integral $I^{\beta}$ of order $\beta\in\R, \beta>0$, by
\begin{equation}
I^{\beta}(f)(t)=\frac1{\Gamma(\beta)}\int_0^t (t-\tau)^{\beta-1}f(\tau)\ d\tau,
\qquad \mbox{a.e.}\ t\in(0,T),
\end{equation}
where $\Gamma (\beta)=\int_0^\infty t^{\beta-1}e^{-t}\ dt$ is the Euler gamma function. 
\end{definition}
We note that
\begin{equation}\label{eq:I1}
I^{1}(f)(t)=\int_0^t f(\tau)\ d\tau\,.
\end{equation}
For the sequel it is convenient to introduce the following function
\begin{equation}
\Phi_{\beta}(t)=\frac{t^{\beta-1}}{\Gamma(\beta)} \quad t>0,
\end{equation}
so
\begin{equation}\label{eq:Iconv}
I^{\beta}(f)(t)=(\Phi_{\beta}*f)(t),
\qquad \mbox{a.e.}\ t\in(0,T).
\end{equation}
For  $f\in L^2(0,T)$ we have
\begin{equation}\label{eq:convol}
\|I^{\beta}(f)\|_{L^2(0,T)}\le \|\Phi_{\beta}\|_{L^1(0,T)}\|f\|_{L^2(0,T)}
\,.
\end{equation}
If we take into account that
\begin{equation}
\Phi_{\beta}*\Phi_{\gamma}(t)=\Phi_{\beta+\gamma}(t) \quad t>0\quad \beta,\gamma>0,
\end{equation}
we have 
\begin{equation}\label{eq:semigroup}
I^{\beta}I^{\gamma}(f)=I^{\beta+\gamma}(f)\,.
\end{equation}

\begin{equation}\label{eq:der-frac}
\partial_t^{\alpha}f(t)=
\begin{cases}
\displaystyle
I^{1-\alpha}\big(\frac{df}{dt}\big)(t)
=\frac1{\Gamma(1-\alpha)}\int_0^t (t-\tau)^{-\alpha}\frac{df}{d\tau}(\tau)\ d\tau\qquad 0<\alpha<1\,,
\\
\\
\displaystyle
I^{2-\alpha}\big(\frac{d^2f}{dt^2}\big)(t)
=\frac1{\Gamma(2-\alpha)}\int_0^t (t-\tau)^{1-\alpha}\frac{d^2f}{d\tau^2}(\tau)\ d\tau\qquad 1<\alpha<2\,.
\end{cases}
\end{equation}
We define the operator $A$ in $L^2(\Omega)$ by
\begin{equation*}
\begin{split}
D(A)&= H^2(\Omega)\cap H^1_0(\Omega)
\\
(Au)(x)&=-\triangle u(x), 
\quad x\in\Omega,
\quad u\in D(A).
\end{split}
\end{equation*}
The fractional powers $A^\theta$ are defined for $\theta>0$, see e.g. \cite{Pazy} and \cite[Example 4.34]{Lunardi}. 
We recall that the spectrum of $A$ consists of a sequence of positive eigenvalues, each of them with finite dimensional eigenspace, and there exists an orthonormal basis of $L^2(\Omega)$ consisting of eigenfunctions of $A$. We denote such a basis by $\{e_n\}_{n\in\N}$ and  by $\lambda_n$ the eigenvalue with eigenfunction $e_n$, that is $Ae_n=\lambda_n e_n$. Then, for $\theta>0$ the domain $D(A^\theta)$ of $A^\theta$ consists of those functions $u\in L^2(\Omega)$ such that
\begin{equation*}
\sum_{n=1}^\infty \lambda_n^{2\theta} |\langle u,e_n\rangle|^2<+\infty
\end{equation*}
and
\begin{equation*}
A^\theta u=\sum_{n=1}^\infty \lambda_n^{\theta} \langle u,e_n\rangle e_n,
\qquad
u\in D(A^\theta).
\end{equation*}
Moreover $D(A^\theta)$ is a Hilbert space with the norm
\begin{equation}\label{eq:norm-frac}
\|u\|_{D(A^\theta)}=\|A^\theta u\|_{L^2(\Omega)}=\left(\sum_{n=1}^\infty \lambda_n^{2\theta} |\langle u,e_n\rangle|^2\right)^{1/2},
\qquad
u\in D(A^\theta)
\,.
\end{equation}
We have $D(A^\theta)\subset H^{2\theta}(\Omega)$. In particular, $D(A^{\frac12})=H^1_0(\Omega)$.
If we identify the dual $(L^{2}(\Omega))'$ with $L^{2}(\Omega)$ itself, then we have $D(A^\theta)\subset L^{2}(\Omega)\subset(D(A^\theta))'$.
From now on we set 
\begin{equation}
D(A^{-\theta}):=(D(A^\theta))',
\end{equation}
whose elements are bounded linear functionals on $D(A^\theta)$. If $\varphi\in D(A^{-\theta})$ and $u\in D(A^\theta)$ the value of $\varphi$ applied to $u$ is denoted  
by 
\begin{equation}\label{eq:duality}
\langle \varphi,u\rangle_{-\theta,\theta}:=\varphi(u)\,.
\end{equation}
 In addition, $D(A^{-\theta})$ is a Hilbert space with the norm
\begin{equation}\label{eq:norm-theta}
\|\varphi\|_{D(A^{-\theta})}=\left(\sum_{n=1}^\infty \lambda_n^{-2\theta} |\langle \varphi,e_n\rangle_{-\theta,\theta}|^2\right)^{1/2},
\qquad
\varphi\in D(A^{-\theta})
\,.
\end{equation}
We also recall that 
\begin{equation}\label{eq:-theta}
\langle \varphi,u\rangle_{-\theta,\theta}=\langle\varphi,u\rangle
\qquad \mbox{for}\ \varphi\in L^{2}(\Omega)\,,u\in D(A^\theta),
\end{equation}
e.g. see \cite[Chapitre V]{B}.

For $\alpha,\beta> 0$ arbitrary constants, we define the Mittag--Leffler functions by
\begin{equation}
E_{\alpha,\beta}(z):= \sum_{k=0}^\infty\frac{z^k}{\Gamma(\alpha k+\beta)},
\quad z\in\C.
\end{equation}
By the power series, one can note that $E_{\alpha,\beta}(z)$ is an entire function
of $z\in\C$.
\begin{lemma}
Let $1<\alpha<2$ and $\beta>0$ be. Then for any $\mu$ such that $\pi\alpha/2<\mu<\pi$ there exists a constant $C=C(\alpha,\beta,\mu)>0$ such that 
\begin{equation}\label{eq:stimeE}
\big|E_{\alpha,\beta}(z)\big|\le \frac{C}{1+|z|},
\qquad \mu\le|\arg(z)|\le\pi.
\end{equation}
\end{lemma}

\begin{lemma} For any $0<\beta<1$ the function $x\to\frac{x^\beta}{1+x}$ gains its maximum on $[0,+\infty[$ at point $\frac\beta{1-\beta}$ and the maximum value is given by
\begin{equation}\label{eq:maxbeta}
\max_{x\ge0}\frac{x^\beta}{1+x}=\beta^\beta(1-\beta)^{1-\beta}
\,.
\end{equation}

\end{lemma}
For a Hilbert space $H$ endowed with the norm $\|\cdot\|_H$ and $\beta\in (0,1)$, $H^\beta(0,T;H)$ is the space of all $u\in L^2(0,T;H)$ such that 
\begin{equation*}
[u]_{H^\beta(0,T;H)}:=\left(\int_0^T\int_0^T\frac{\|u(t)-u(\tau)\|_H^2}{|t-\tau|^{1+2\beta}}\ dtd\tau\right)^{1/2}<+\infty\,,
\end{equation*}
that is $[u]_{H^\beta(0,T;H)}$ is the so-called Gagliardo semi-norm of $u$.
$H^\beta(0,T;H)$ is endowed with the norm
\begin{equation}\label{eq:defHs}
\|\cdot\|_{H^\beta(0,T;H)}:=\|\cdot\|_{L^2(0,T;H)}+[\ \cdot\ ]_{H^\beta(0,T;H)}.
\end{equation}
We will use later the following extension to the case of vector valued functions 
of a known result, see
\cite[Theorem 2.1]{GLY}.
\begin{theorem}\label{th:R-Lop}
Let $H$ be a separable Hilbert space.
\begin{itemize}
\item[(i)] 
The Riemann--Liouville operator $I^{\beta}:L^2(0,T;H)\to L^2(0,T;H)$, $0<\beta\le1$, is injective and the range ${\cal R}(I^{\beta})$ of $I^{\beta}$ is given by
\begin{equation}\label{eq:range}
{\cal R}(I^{\beta})=
\begin{cases}H^\beta(0,T;H), \hskip5.5cm 0<\beta<\frac12,
\\
\Big\{v\in H^{\frac12}(0,T;H): \int_0^Tt^{-1}|v(t)|^2 dt<\infty\Big\},
\qquad \beta=\frac12,
\\
_0H^\beta(0,T;H), \hskip5.2cm \frac12<\beta\le1,
\end{cases}
\end{equation}
where $_0H^\beta(0,T)=\{u\in H^\beta(0,T): u(0)=0\}$.
\item[(ii)] 
For the Riemann--Liouville operator $I^{\beta}$ and its inverse operator $I^{-\beta}$ the norm equivalences 
\begin{equation}\label{eq:R-Lop}
\begin{split}
\|I^{\beta}(u)\|_{H^\beta(0,T;H)}
&\sim\|u\|_{L^2(0,T;H)},
\qquad u\in L^2(0,T;H),
\\
\|I^{-\beta}(v)\|_{L^2(0,T;H)}
&\sim\|v\|_{H^\beta(0,T;H)},
\qquad v\in {\cal R}(I^{\beta}),
\end{split}
\end{equation}
hold true.
\end{itemize}
\end{theorem}
For the sake of completeness, we recall the notion of a weak solution for fractional diffusion-wave equations,  see \cite[Definition 2.1]{SY}.
\begin{definition}
Let $1<\alpha<2$.
We define $u$ as a weak solution to problem
\begin{equation}\label{eq:weakp}
\begin{cases}
\displaystyle
\partial_t^{\alpha}u(t,x)=\triangle u(t,x)
\qquad t\in (0,T),\  x\in\Omega,
\\
u(t,x)=0    \hskip2.6cm t\in (0,T),\ x\in\partial\Omega,
\\
u(0,x)=u_0(x), \quad u_t(0,x)=u_1(x),\quad x\in\Omega,
\end{cases}
\end{equation}
if  $\partial_t^{\alpha}u(t,\cdot)=\triangle u(t,\cdot)$ holds in $L^2(\Omega)$, $u(t,\cdot)\in H_0^1(\Omega)$ for almost all $t\in (0,T)$ and for some $\theta>0$, depending on the initial data $u_0,u_1$, one has
$u, \partial_t u\in C([0,T];D(A^{-\theta}))$  and
\begin{equation}
\lim_{t\to0}\|u(t,\cdot)-u_0\|_{D(A^{-\theta})}=\lim_{t\to0}\|\partial_tu(t,\cdot)-u_1\|_{D(A^{-\theta})}=0\,.
\end{equation}
\end{definition}
We also need to recall some existence results  given in \cite[Theorem 2.3]{SY}, that we have integrated with other essential regularity properties of the solution, see \eqref{eq:error} below.
\begin{theorem}\label{th:saka-y}
\begin{itemize}
\item[(i)]
Let $u_0\in L^2(\Omega)$ and $u_1\in D(A^{-\frac1\alpha})$. Then there exists a unique weak solution $u\in C([0,T];L^2(\Omega))\cap C((0,T];H^2(\Omega)\cap H^1_0(\Omega))$ to \eqref{eq:weakp} with $\partial_t^{\alpha}u\in C((0,T];L^2(\Omega))$ and satisfying 
\begin{equation*}
\lim_{t\to0}\|u(t,\cdot)-u_0\|_{L^2(\Omega)}=0\,,
\quad
\|u\|_{C([0,T];L^2(\Omega))}
\le C\big(\|u_0\|_{L^2(\Omega)}+\|u_1\|_{D(A^{-\frac1\alpha})}\big)
\,,
\end{equation*}
\begin{equation}\label{eq:error}
\begin{split}
&\lim_{t\to0}\|\partial_tu(t,\cdot)-u_1\|_{D(A^{-\theta})}=0\,,
\qquad \theta\in\Big(\frac1\alpha,1\Big),
\\
&\|\partial_tu\|_{C([0,T];D(A^{-\theta}))}
\le C\big(\|u_0\|_{L^2(\Omega)}+\|u_1\|_{D(A^{-\frac1\alpha})}\big)
\,,
\end{split}
\end{equation}
for some constant $C>0$.
Moreover, if $u_1\in L^2(\Omega)$ we have
\begin{align}
\label{eq:def-u0}
u(t,x)&=\sum_{n=1}^\infty\big[ \langle u_0,e_n\rangle E_{\alpha,1}(-\lambda_nt^\alpha)
+\langle u_1,e_n\rangle t E_{\alpha,2}(-\lambda_nt^\alpha)\big]e_n(x),
\\
\label{eq:def-u_t}
\partial_tu(t,x)&=\sum_{n=1}^\infty\big[-\lambda_n \langle u_0,e_n\rangle t^{\alpha-1}E_{\alpha,\alpha}(-\lambda_nt^\alpha)
+\langle u_1,e_n\rangle E_{\alpha,1}(-\lambda_nt^\alpha)\big]e_n(x),
\\
\label{eq:def-u-alpha}
\partial_t^{\alpha}u(t,x)&=\sum_{n=1}^\infty\big[-\lambda_n \langle u_0,e_n\rangle E_{\alpha,1}(-\lambda_nt^\alpha)
-\lambda_n\langle u_1,e_n\rangle t E_{\alpha,2}(-\lambda_nt^\alpha)\big]e_n(x)
\,,
\end{align}
\begin{equation*}
\|\partial_tu(t,\cdot)\|_{L^2(\Omega)}
\le C\big(t^{-1}\|u_0\|_{L^2(\Omega)}+\|u_1\|_{L^2(\Omega)}\big)
\qquad (C>0)\,.
\end{equation*}
\item[(ii)]
If $u_0\in H^2(\Omega)\cap H^1_0(\Omega)$ and $u_1\in H^1_0(\Omega)$, then the unique weak solution $u$ to \eqref{eq:weakp} given by \eqref{eq:def-u0} belongs to $C([0,T];H^2(\Omega)\cap H^1_0(\Omega))\cap C^1([0,T];L^2(\Omega))$ and $\partial_t^{\alpha}u\in C([0,T];L^2(\Omega))$. In addition, there exists a constant $C>0$ such that
\begin{equation}
\|u\|_{C([0,T];H^2(\Omega))}+\|u\|_{C^1([0,T];L^2(\Omega))}+\|\partial_t^{\alpha}u\|_{C([0,T];L^2(\Omega))}
\le C\big(\|u_0\|_{H^2(\Omega)}+\|u_1\|_{H^1(\Omega)}\big)
\,.
\end{equation}

\end{itemize}
\end{theorem}
\begin{Proof}
We refer to  \cite[Theorem 2.3]{SY} for the proof of all statements, except for the proof of \eqref{eq:error}.
We first observe that,
since $u_1\in D(A^{-\frac1\alpha})$, the expression \eqref{eq:def-u_t} for $\partial_tu$ has to be written in the form
\begin{equation*}
\partial_tu(t,x)=\sum_{n=1}^\infty\big[-\lambda_n \langle u_0,e_n\rangle t^{\alpha-1}E_{\alpha,\alpha}(-\lambda_nt^\alpha)
+\langle u_1,e_n\rangle_{-\frac1\alpha,\frac1\alpha} E_{\alpha,1}(-\lambda_nt^\alpha)\big]e_n(x)\,.
\end{equation*}
For $\theta\in(0,1)$ to choose suitably later, we have
\begin{multline}\label{eq:error1}
\|\partial_tu(t,\cdot)-u_1\|_{D(A^{-\theta})}^2
\\
=
\sum_{n=1}^\infty\lambda_n^{-2\theta}\big|-\lambda_n \langle u_0,e_n\rangle t^{\alpha-1}E_{\alpha,\alpha}(-\lambda_nt^\alpha)
+\langle u_1,e_n\rangle_{-\frac1\alpha,\frac1\alpha} \big(E_{\alpha,1}(-\lambda_nt^\alpha)-1\big)\big|^2
\\
\le
2t^{2(\alpha-1)}\sum_{n=1}^\infty\lambda_n^{2(1-\theta)}\ | \langle u_0,e_n\rangle E_{\alpha,\alpha}(-\lambda_nt^\alpha)|^2
+2\sum_{n=1}^\infty\lambda_n^{-2\theta}\big|\langle u_1,e_n\rangle_{-\frac1\alpha,\frac1\alpha} \big(E_{\alpha,1}(-\lambda_nt^\alpha)-1\big)\big|^2
\,.\end{multline}
To estimate the first sum we use  \eqref{eq:stimeE} and \eqref{eq:maxbeta} to get
\begin{equation*}
t^{2(\alpha-1)}\lambda_n^{2(1-\theta)}\ | \langle u_0,e_n\rangle E_{\alpha,\alpha}(-\lambda_nt^\alpha)|^2
\le C
t^{2(\alpha\theta-1)}\Big(\frac{(\lambda_nt^{\alpha})^{1-\theta}}{1+\lambda_nt^\alpha}\Big)^2 | \langle u_0,e_n\rangle|^2
\le C
t^{2(\alpha\theta-1)} | \langle u_0,e_n\rangle|^2
\,,
\end{equation*}
while, regarding the second sum, we have
\begin{equation*}
\lambda_n^{-2\theta}\big|\langle u_1,e_n\rangle_{-\frac1\alpha,\frac1\alpha} \big(E_{\alpha,1}(-\lambda_nt^\alpha)-1\big)\big|^2
=
\lambda_n^{-2(\theta-\frac1\alpha)}\lambda_n^{-\frac2\alpha}|\langle u_1,e_n\rangle_{-\frac1\alpha,\frac1\alpha}|^2 \big|E_{\alpha,1}(-\lambda_nt^\alpha)-1\big|^2
\,.
\end{equation*}
Therefore, plugging the above two estimates into \eqref{eq:error1} we obtain
\begin{multline*}
\|\partial_tu(t,\cdot)-u_1\|_{D(A^{-\theta})}^2
\\
\le
Ct^{2(\alpha\theta-1)}\|u_0\|_{L^2(\Omega)}^2
+2\sum_{n=1}^\infty\lambda_n^{-2(\theta-\frac1\alpha)}\lambda_n^{-\frac2\alpha}|\langle u_1,e_n\rangle_{-\frac1\alpha,\frac1\alpha}|^2 \big|E_{\alpha,1}(-\lambda_nt^\alpha)-1\big|^2
\,,
\end{multline*}
whence it follows that for $\theta>\frac1\alpha$ \eqref{eq:error} holds true.
\end{Proof}

%

\section{Regularity for $u_0\in H^1_0(\Omega)$ and $u_1\in L^2(\Omega)$}
We establish a result about the regularity of the weak solutions assuming on the data $u_0,u_1$  a degree of regularity intermediate between those assumed in (i) and (ii) of Theorem \ref{th:saka-y}.
\begin{theorem}\label{th:reg-l2}
If $u_0\in H^1_0(\Omega)$ and $u_1\in L^2(\Omega)$, then the unique weak solution $u$ to \eqref{eq:weakp} given by \eqref{eq:def-u0}--\eqref{eq:def-u-alpha} belongs to $C([0,T];H^1_0(\Omega))\cap C^1([0,T];D(A^{-\theta}))$, $\theta\in\big(\frac{2-\alpha}{2\alpha},\frac12\big]$, 
and
\begin{equation}\label{eq:in-data}
\begin{split}
&\lim_{t\to0}\|u(t,\cdot)-u_0\|_{H^1_0(\Omega)}=
\lim_{t\to0}\|\partial_t u(t,\cdot)-u_1\|_{D(A^{-\theta})}=0\,,
\\
&\|u\|_{C([0,T];H^1_0(\Omega))}+\|\partial_t u\|_{C([0,T];D(A^{-\theta}))}
\le
C(\| u_0\|_{H^1_0(\Omega)}+\| u_1\|_{L^2(\Omega)}).
\end{split}
\end{equation}
In addition, for any $\theta\in\big(0,\frac1{2\alpha}\big)$ there exists a constant $C>0$ such that
\begin{equation}\label{eq:nabla-l2}
\|\nabla u\|_{L^2(0,T;D(A^{\theta}))}\le C\big(\|u_0\|_{H^1_0(\Omega)}+\|u_1\|_{L^2(\Omega)}\big),
\end{equation}
and for any $\theta\in\big(\frac{\alpha-1}{2\alpha},\frac12\big)$ there exists a constant $C>0$ such that
\begin{equation}\label{eq:partial-l2}
\|\partial_t^{\alpha}u\|_{L^2(0,T;D(A^{-\theta}))}\le C\big(\|u_0\|_{H^1_0(\Omega)}+\|u_1\|_{L^2(\Omega)}\big)\,.
\end{equation}
Moreover, if we assume $u_0\in D(A^{\frac12+\varepsilon})$ with $\varepsilon\in\big(\frac{2-\alpha}{2\alpha},\frac12\big)$, then
\begin{equation}\label{eq:in-data1}
\begin{split}
&
\lim_{t\to0}\|\partial_t u(t,\cdot)-u_1\|_{L^2(\Omega)}=0\,,
\\
&\|\partial_t u\|_{C([0,T];L^2(\Omega))}
\le
C(\| u_0\|_{D(A^{\frac12+\varepsilon})}+\| u_1\|_{L^2(\Omega)}).
\end{split}
\end{equation}
\end{theorem}
\begin{Proof}
In virtue of the expression \eqref{eq:def-u0} for the solution $u$
we have
\begin{multline}\label{eq:u000}
\|u(t,\cdot)-u_0\|_{H^1_0(\Omega)}^2
=
\sum_{n=1}^\infty\lambda_n
\big| \langle u_0,e_n\rangle \big(E_{\alpha,1}(-\lambda_nt^\alpha)-1\big)+\langle u_1,e_n\rangle t E_{\alpha,2}(-\lambda_nt^\alpha)\big|^2\\
\le
2\sum_{n=1}^\infty\lambda_n
\big| \langle u_0,e_n\rangle\big|^2 \big|E_{\alpha,1}(-\lambda_nt^\alpha)-1\big|^2
+t^{2-\alpha} 2C^2\sum_{n=1}^\infty\big| \langle u_1,e_n\rangle\big|^2 
\Big(\frac{(\lambda_nt^{\alpha})^{\frac12}}{1+\lambda_nt^\alpha}\Big)^2,
\end{multline}
thanks also to \eqref{eq:stimeE}.
We observe that for any $n\in\N$ $\lim_{t\to0}\big(E_{\alpha,1}(-\lambda_nt^\alpha)-1\big)=0$. Moreover,
again by \eqref{eq:stimeE}, we get for $n\in\N$ and $0\le t\le T$
\begin{equation*}
\lambda_n
\big| \langle u_0,e_n\rangle\big|^2 \big|E_{\alpha,1}(-\lambda_nt^\alpha)-1\big|^2
\le
2\lambda_n
\big| \langle u_0,e_n\rangle\big|^2 \Big(\frac{C}{(1+\lambda_nt^\alpha)^2}+1\Big)
\le
C\lambda_n\big| \langle u_0,e_n\rangle\big|^2,
\end{equation*}
hence by \eqref{eq:u000} we deduce $\lim_{t\to0}\|u(t,\cdot)-u_0\|_{H^1_0(\Omega)}=0$ and for any $t\in[0,T]$
\begin{equation*}
\|u(t,\cdot)\|_{H^1_0(\Omega)}^2
\le
C(\| u_0\|_{H^1_0(\Omega)}^2+\| u_1\|_{L^2(\Omega)}^2).
\end{equation*}
To complete the proof of \eqref{eq:in-data}, we fix $\theta\in\big(\frac{2-\alpha}{2\alpha},\frac12\big]$ and use formula \eqref{eq:def-u_t} to note that 
\begin{multline}\label{eq:reg.u-t}
\|\partial_t u(t,\cdot)-u_1\|_{D(A^{-\theta})}^2=
\sum_{n=1}^\infty\lambda_n^{-2\theta}\big|-\lambda_n \langle u_0,e_n\rangle t^{\alpha-1}E_{\alpha,\alpha}(-\lambda_nt^\alpha)
+\langle u_1,e_n\rangle \big(E_{\alpha,1}(-\lambda_nt^\alpha)-1\big)\big|^2
\\
\le
Ct^{\alpha-2+2\alpha\theta}\sum_{n=1}^\infty\lambda_n\big| \langle u_0,e_n\rangle\big|^2 \left(\frac{(\lambda_nt^{\alpha})^{\frac{1-2\theta}2}}{1+\lambda_nt^\alpha}\right)^2
+2\sum_{n=1}^\infty\lambda_n^{-2\theta}\big|\langle u_1,e_n\rangle\big|^2 \big|E_{\alpha,1}(-\lambda_nt^\alpha)-1\big|^2,
\end{multline}
thanks also to \eqref{eq:stimeE}.
Since $0<\frac{1-2\theta}2<1$ we can apply \eqref{eq:maxbeta} to have
\begin{equation*}
\|\partial_t u(t,\cdot)-u_1\|_{D(A^{-\theta})}^2
\le
Ct^{\alpha-2+2\alpha\theta}\| u_0\|_{H^1_0(\Omega)}^2
+2\sum_{n=1}^\infty\big|\langle u_1,e_n\rangle\big|^2 \big|E_{\alpha,1}(-\lambda_nt^\alpha)-1\big|^2.
\end{equation*}
Therefore, by analogous argumentations to those done before, since $\alpha-2+2\alpha\theta>0$ we deduce 
\break
$\lim_{t\to0}\|\partial_t u(t,\cdot)-u_1\|_{D(A^{-\theta})}=0$ and for any $t\in[0,T]$
\begin{equation*}
\|\partial_t u(t,\cdot)\|_{D(A^{-\theta})}^2
\le
C(\| u_0\|_{H^1_0(\Omega)}^2+\| u_1\|_{L^2(\Omega)}^2).
\end{equation*}

\begin{multline*}
\|\nabla u(\cdot,t)\|_{D(A^{\theta})}^2=\sum_{n=1}^\infty\lambda_n^{1+2\theta}\big| \langle u_0,e_n\rangle E_{\alpha,1}(-\lambda_nt^\alpha)
+\langle u_1,e_n\rangle t E_{\alpha,2}(-\lambda_nt^\alpha)\big|^2
\\
\le
C\sum_{n=1}^\infty\lambda_n |\langle u_0,e_n\rangle|^2 \frac{\lambda_n^{2\theta}}{(1+\lambda_nt^\alpha)^2}
+C\sum_{n=1}^\infty|\langle u_1,e_n\rangle|^2 \frac{\lambda_n^{1+2\theta}t^2}{(1+\lambda_nt^\alpha)^2}
\,.
\end{multline*}
Since
\begin{equation*}
 \frac{\lambda_n^{2\theta}}{(1+\lambda_nt^\alpha)^2}=\Big(\frac{(\lambda_nt^\alpha)^{\theta}}{1+\lambda_nt^\alpha}\Big)^2 t^{-2\alpha\theta},
\qquad
\frac{\lambda_n^{1+2\theta}t^2}{(1+\lambda_nt^\alpha)^2}=\Big(\frac{(\lambda_nt^\alpha)^{\frac{1+2\theta}2}}{1+\lambda_nt^\alpha}\Big)^2 t^{2-\alpha(1+2\theta)},
\end{equation*}
for $0<\theta<\frac12$, we  can apply \eqref{eq:maxbeta} to have
\begin{equation}
\|\nabla u(\cdot,t)\|_{L^2(0,T;D(A^{\theta}))}^2
\le
Ct^{-2\alpha\theta}\|u_0\|_{H^1_0(\Omega)}^2
+Ct^{2-\alpha(1+2\theta)}\|u_1\|_{L^2(\Omega)}^2
\end{equation}
Thanks to \eqref{eq:norm-theta}, \eqref{eq:def-u-alpha} and \eqref{eq:stimeE} we get
\begin{multline}\label{eq:}
\|\partial_t^{\alpha}u(\cdot,t)\|_{D(A^{-\theta})}^2
=\sum_{n=1}^\infty\lambda_n^{-2\theta}\big|\lambda_n \langle u_0,e_n\rangle E_{\alpha,1}(-\lambda_nt^\alpha)
+\lambda_n\langle u_1,e_n\rangle t E_{\alpha,2}(-\lambda_nt^\alpha)\big|^2
\\
\le
C\sum_{n=1}^\infty\lambda_n |\langle u_0,e_n\rangle|^2 \frac{\lambda_n^{1-2\theta}}{(1+\lambda_nt^\alpha)^2}
+C\sum_{n=1}^\infty|\langle u_1,e_n\rangle|^2 \frac{\lambda_n^{2(1-\theta)}t^2}{(1+\lambda_nt^\alpha)^2}
\end{multline}
\begin{equation*}
\frac{\lambda_n^{1-2\theta}}{(1+\lambda_nt^\alpha)^2}=\Big(\frac{(\lambda_nt^\alpha)^{\frac{1-2\theta}2}}{1+\lambda_nt^\alpha}\Big)^2 t^{\alpha(2\theta-1)}
\end{equation*}
\begin{equation*}
\frac{\lambda_n^{2(1-\theta)}t^2}{(1+\lambda_nt^\alpha)^2}
=\Big(\frac{(\lambda_nt^\alpha)^{1-\theta}}{1+\lambda_nt^\alpha}\Big)^2 t^{2+2\alpha(\theta-1)}
\end{equation*}

\begin{equation}\label{eq:}
\|\partial_t^{\alpha}u(\cdot,t)\|_{D(A^{-\frac1{2\alpha}})}^2
\le C t^{1-\alpha}\ \| u_0\|_{H^1_0(\Omega)}^2+Ct ^{3-2\alpha}\| u_1\|_{L^2(\Omega)}^2\,.
\end{equation}

By assuming, in addition, that $u_0\in D(A^{\frac12+\varepsilon})$ with $\varepsilon\in\big(\frac{2-\alpha}{2\alpha},\frac12\big)$ we have
\begin{multline}
\|\partial_t u(t,\cdot)-u_1\|_{L^2(\Omega)}^2=
\sum_{n=1}^\infty\big|-\lambda_n \langle u_0,e_n\rangle t^{\alpha-1}E_{\alpha,\alpha}(-\lambda_nt^\alpha)
+\langle u_1,e_n\rangle \big(E_{\alpha,1}(-\lambda_nt^\alpha)-1\big)\big|^2
\\
\le
Ct^{\alpha-2+2\alpha\varepsilon}\sum_{n=1}^\infty\lambda_n^{1+2\varepsilon}\big| \langle u_0,e_n\rangle\big|^2 \left(\frac{(\lambda_nt^{\alpha})^{\frac{1-2\varepsilon}2}}{1+\lambda_nt^\alpha}\right)^2
+2\sum_{n=1}^\infty\big|\langle u_1,e_n\rangle\big|^2 \big|E_{\alpha,1}(-\lambda_nt^\alpha)-1\big|^2.
\end{multline}
Thanks to \eqref{eq:maxbeta} with $\beta=\frac{1-2\varepsilon}2$ we obtain
\begin{equation*}
\|\partial_t u(t,\cdot)-u_1\|_{L^2(\Omega)}^2
\le
Ct^{\alpha-2+2\alpha\varepsilon}\| u_0\|_{D(A^{\frac12+\varepsilon})}
+2\sum_{n=1}^\infty\big|\langle u_1,e_n\rangle\big|^2 \big|E_{\alpha,1}(-\lambda_nt^\alpha)-1\big|^2,
\end{equation*}
hence, since $\alpha-2+2\alpha\varepsilon>0$, we deduce \eqref{eq:in-data1}.
\end{Proof}

\begin{remark}
Comparing the regularity results given in Theorems \ref{th:saka-y} and \ref{th:reg-l2},
we have to observe that if $\theta\in\big(\frac{2-\alpha}{2\alpha},\frac12\big]$ then $D(A^{-\theta})\subset D(A^{-\eta})$ for any 
$\eta\in\big(\frac1{\alpha},1\big]$. Therefore Theorem \ref{th:reg-l2}  effectively improves the regularity of the weak solution. 

Moreover, taking into account the argumentations used to get \eqref{eq:reg.u-t}, we note that to secure a regularity of $\partial_t u$ in $L^2(\Omega)$ we have to assume the datum $u_0$ more regular than $u_0\in H^1_0(\Omega)=D(A^{\frac12})$, that is $u_0\in D(A^{\frac12+\varepsilon})$ with $\varepsilon\in\big(\frac{2-\alpha}{2\alpha},\frac12\big)$, see \eqref{eq:in-data1}.

\end{remark}

\section{Hidden regularity results}
To begin with we single out some technical results that we will use later in the main theorem.
\begin{lemma}\label{le:tech0}
For any $w\in H^2(\Omega)$ one has
\begin{multline}\label{eq:triangleuF00}
2\int_\Omega
\triangle w\
h\cdot\nabla w\ dx 
=
\int_{\partial\Omega}\Big[2\partial_\nu w\ h\cdot\nabla w-h\cdot\nu |\nabla w|^2\Big]\ d\sigma 
-2
\sum_{i,j=1}^N\int_\Omega
\partial_i  h_j\partial_i w\partial_j w\ dx 
\\
\int_{\Omega}
\sum_{j=1}^N \partial_jh_j\ |\nabla w|^2\ dx 
\,.
\end{multline}
\end{lemma}
\begin{Proof}
We integrate by parts to get
\begin{equation}\label{eq:triangleu00}
\int_\Omega
\triangle w\
h\cdot\nabla w\ dx
=
\int_{\partial\Omega}
\partial_\nu w\
h\cdot\nabla w\ d\sigma
-\int_\Omega\nabla w
\cdot\nabla \big( h\cdot\nabla w\big)\ dx
\,.
\end{equation}
Since
\begin{multline*}
\int_\Omega\nabla w\cdot\nabla \big( h\cdot\nabla w \big)\ dx
\\
=\sum_{i,j=1}^N\int_\Omega
\partial_i w\ \partial_i ( h_j\partial_j w)\ dx
=\sum_{i,j=1}^N\int_\Omega\partial_i w\ \partial_i  h_j\partial_j w\ dx
+\sum_{i,j=1}^N\int_\Omega h_j\ \partial_i w\partial_j ( \partial_i w)\ dx,
\end{multline*}
we evaluate the last term on the right-hand side again by an integration  by parts, so we obtain
\begin{equation*}
\begin{split}
\sum_{i,j=1}^N\int_\Omega
h_j\ \partial_i w \partial_j ( \partial_i w)\ dx
=&\frac12\sum_{j=1}^N\int_\Omega
h_j\ \partial_j \Big( \sum_{i=1}^N(\partial_i w)^2\Big)\ dx
\\
=&\frac12\int_{\partial\Omega}
h\cdot\nu |\nabla w|^2\ d\sigma
-\frac12\int_{\Omega}
\sum_{j=1}^N \partial_jh_j\ |\nabla w|^2\ dx
\,.
\end{split}
\end{equation*}
Therefore, if we merge the above two identities with \eqref{eq:triangleu00}, then we have \eqref{eq:triangleuF00}.
\end{Proof}
\begin{lemma}\label{le:tech}
Assume $1<\alpha<2$ and the weak solution $u$ of 
\begin{equation}\label{eq:stato}
\partial_t^{\alpha}u(t,x) =\triangle u (t,x)
\quad 
\text{in}
\ \ 
(0,\infty)\times\Omega
\end{equation}
belonging to $C([0,+\infty);H^2(\Omega)\cap H^1_0(\Omega))\cap C^1([0,+\infty);L^2(\Omega))$ with $\partial_t^{\alpha}u\in C([0,+\infty);L^2(\Omega))$.
Then, for a vector field $h:\overline{\Omega}\to\R^N$ of class $C^1$ and $\beta,\theta\in(0,1)$ the following identities hold true
\begin{multline}\label{eq:identity}
\int_{\partial\Omega}
\Big[2I^{\beta}(\partial_\nu u)(t)\  h\cdot I^{\beta}(\nabla u)(t)
-h\cdot\nu \big|  I^{\beta}(\nabla u)(t)\big|^2\Big]
\ d\sigma
=2\langle I^{\beta}(\partial_t^{\alpha}u)(t),h\cdot I^{\beta}(\nabla u)(t)\rangle_{-\theta,\theta}
\\
+2\sum_{i,j=1}^N\int_\Omega
\partial_i  h_j I^{\beta}(\partial_i u)(t)I^{\beta}(\partial_j u)(t)\ dx
-\int_{\Omega}
\sum_{j=1}^N \partial_jh_j\ | I^{\beta}(\nabla u)(t)|^2\ dx
\,, \qquad t>0,
\end{multline}
\begin{multline}\label{eq:identity2}
\int_{\partial\Omega}\Big[2 \big(I^{\beta}(\partial_\nu u )(t)-I^{\beta}(\partial_\nu u )(\tau)\big)\ h\cdot\big(I^{\beta}(\nabla u)(t)-I^{\beta}(\nabla u)(\tau)\big)
-h\cdot\nu \big| I^{\beta}(\nabla u )(t)-I^{\beta}(\nabla u )(\tau)\big|^2\Big] d\sigma 
\\
=
2\langle I^{\beta}(\partial_t^{\alpha}u)(t)-I^{\beta}(\partial_{t}^{\alpha}u)(\tau),h\cdot \big(I^{\beta}(\nabla u)(t)-I^{\beta}(\nabla u)(\tau)\big)\rangle_{-\theta,\theta}
\\
+2
\sum_{i,j=1}^N\int_\Omega
\partial_i  h_j  \big(I^{\beta}(\partial_i u)(t)-I^{\beta}(\partial_i u)(\tau)\big)  \big(I^{\beta}(\partial_j u)(t)-I^{\beta}(\partial_j u)(\tau)\big)\ dx 
\\
-
\sum_{j=1}^N\int_{\Omega} \partial_jh_j\ |I^{\beta}(\nabla u )(t)-I^{\beta}(\nabla u )(\tau)|^2\ dx 
\,, \qquad t,\tau>0
\,.
\end{multline}
\end{lemma}
\begin{Proof}
First, we apply the operator $I^{\beta}$, $\beta\in(0,1)$, to  equation \eqref{eq:stato}:
\begin{equation}\label{eq:eqI}
I^{\beta}(\partial_t^{\alpha}u)(t)=I^{\beta}(\triangle u )(t)
\qquad t>0.
\end{equation}
Fix $\theta\in(0,1)$, 
by means of the duality 
$\langle \cdot,\cdot\rangle_{-\theta,\theta}$ introduced by \eqref{eq:duality}
we multiply the terms of the previous equation by 
\begin{equation*}
2h\cdot \nabla I^{\beta}( u)(t),
\end{equation*}
that is
\begin{equation*}
2\langle I^{\beta}(\partial_t^{\alpha}u)(t),h\cdot\nabla I^{\beta}(u)(t)\rangle_{-\theta,\theta}
=2\langle \triangle I^{\beta}( u )(t),h\cdot\nabla I^{\beta}(u)(t)\rangle_{-\theta,\theta}
.
\end{equation*}
Thanks to the regularity of data and \eqref{eq:-theta} the term on the right-hand side of the previuos equation can be written as a scalar product in $L^2(\Omega)$, so we have
\begin{equation}\label{eq:identity1}
2\langle I^{\beta}(\partial_t^{\alpha}u)(t),h\cdot\nabla I^{\beta}(u)(t)\rangle_{-\theta,\theta}
=
2\int_\Omega \triangle I^{\beta}( u )(t)
h\cdot\nabla I^{\beta}(u)(t)\ dx  
\end{equation}
To evaluate the term
\begin{equation*}
2\int_\Omega\triangle I^{\beta}( u )(t)
h\cdot\nabla I^{\beta}(u)(t)\ dx  
\,,
\end{equation*} 
we apply Lemma \ref{le:tech0} to the function
$
w(t,x)=I^{\beta}( u )(t)
$,
so
from \eqref{eq:triangleuF00} we deduce
\begin{multline*}
2\int_\Omega
\triangle I^{\beta}( u )(t)
h\cdot\nabla I^{\beta}(u)(t)\ dx  
=\int_{\partial\Omega}\Big[2 I^{\beta}(\partial_\nu u )(t)\ h\cdot I^{\beta}(\nabla u)(t)
-h\cdot\nu \big| I^{\beta}(\nabla u )(t)\big|^2\Big] d\sigma 
\\
-2
\sum_{i,j=1}^N\int_\Omega
\partial_i  h_j  I^{\beta}(\partial_i u)(t)  I^{\beta}(\partial_j u)(t)\ dx 
+
\int_{\Omega} \sum_{j=1}^N\partial_jh_j\ |I^{\beta}(\nabla u )(t)|^2\ dx 
\,.
\end{multline*}
In conclusion, plugging the above formula into \eqref{eq:identity1}, we obtain \eqref{eq:identity}.

The proof of \eqref{eq:identity2} is similar to that of \eqref{eq:identity}. Indeed, starting from 
\begin{equation*}
I^{\beta}(\partial_t^{\alpha}u)(t)-I^{\beta}(\partial_{t}^{\alpha}u)(\tau)=I^{\beta}(\triangle u )(t)-I^{\beta}(\triangle u )(\tau)
\qquad t,\tau>0,
\end{equation*}
by means of the duality $\langle \cdot,\cdot\rangle_{-\theta,\theta}$ one multiplies both terms   
by 
\begin{equation*}
2h\cdot\nabla \big(I^{\beta}(u)(t)-I^{\beta}(u)(\tau)\big)\,.
\end{equation*}
Then applying  Lemma \ref{le:tech0} to the function
$
w(t,\tau,x)=I^{\beta}( u )(t)-I^{\beta}( u )(\tau)
$, one can get the  identity \eqref{eq:identity2}. We omit the details.
\end{Proof}
\begin{theorem}
Let  $u_0\in H^2(\Omega)\cap H^1_0(\Omega)$, $u_1\in H^1_0(\Omega)$ and $u$ the weak solution of 
\begin{equation}\label{eq:cauchy1}
\begin{cases}
\displaystyle
\partial_t^{\alpha}u(t,x) =\triangle u (t,x)\,,
\quad t\ge0,\,\, x\in \Omega,
\\
u(t,x)=0\qquad   t\ge0, \,\, x\in\partial\Omega,
\\
u(0,x)=u_{0}(x),\quad
u_t(0,x)=u_{1}(x),\qquad  x\in \Omega.
\end{cases}
\end{equation}
Then, for any $T>0$ there is a constant $c_0=c_0(T)$ such that $u$ satisfies the inequality
\begin{equation}\label{eq:hidden-alpha}
\int_0^T\int_{\partial\Omega} \big|\partial_\nu u\big|^2d\sigma dt
\le c_0(\|\nabla u_0\|^2_{L^2(\Omega)}+\|u_1\|^2_{L^2(\Omega)})
\,.
\end{equation}

\end{theorem}

\begin{Proof}
We will use Theorem \ref{th:R-Lop} with $H=L^2(\partial\Omega)$ and $\beta\in(0,1)$. Indeed,
thanks to \eqref{eq:R-Lop}  we have
\begin{equation}
\|\partial_\nu u\|_{L^2(0,T;L^2(\partial\Omega))}\sim\|I^{\beta}(\partial_\nu u)\|_{H^{\beta}(0,T;L^2(\partial\Omega))}\,,
\end{equation}
so, taking also into account \eqref{eq:defHs}, the proof of \eqref{eq:hidden-alpha} is equivalent to prove
\begin{equation}\label{eq:hidden-alpha1}
\big\|I^{\beta}(\partial_\nu u)\big\|_{L^2(0,T;L^2(\partial\Omega))}^2
+\big[ I^{\beta}(\partial_\nu u) \big]_{H^{\beta}(0,T;L^2(\partial\Omega))}^2
\le c_0(\|\nabla u_0\|^2_{L^2(\Omega)}+\|u_1\|^2_{L^2(\Omega)})
\,.
\end{equation}
To this end we will employ the two identities in Lemma \ref{le:tech} with a suitable choice of the vector field $h$. 
Indeed, we take a vector field $h\in C^1(\overline{\Omega};\R^N)$ satisfying the condition
\begin{equation}\label{eq:h}
h=\nu
\qquad
\text{on}\quad\partial\Omega
\end{equation}
(see e.g. \cite {K} for the existence of such vector field $h$) and first consider the identity \eqref{eq:identity}.
Since
\begin{equation}\label{eq:mmm}
\nabla u=(\partial_\nu u)\nu
\quad
\text{on}
\quad (0,T)\times\partial\Omega\,,
\end{equation}
(see e.g. \cite[Lemma 2.1]{MM} for a detailed proof)
the left-hand side of \eqref{eq:identity} becomes
\begin{equation*}
\int_{\partial\Omega} \big|I^{\beta}(\partial_\nu u)\big|^2d\sigma 
\,.
\end{equation*}
Thanks to that choice of $h$, if we
integrate  \eqref{eq:identity} over $[0,T]$, then we obtain
\begin{multline}\label{eq:identity01}
\int_0^T\int_{\partial\Omega}
\big|I^{\beta}(\partial_\nu u)\big|^2
\ d\sigma dt
=2\int_0^T\langle I^{\beta}(\partial_t^{\alpha}u)(t),h\cdot I^{\beta}(\nabla u)(t)\rangle_{-\theta,\theta}\ dt
\\
+2\sum_{i,j=1}^N\int_0^T\int_\Omega
\partial_i  h_j I^{\beta}(\partial_i u)(t)I^{\beta}(\partial_j u)(t)\ dx dt
-\int_0^T\int_{\Omega}
\sum_{j=1}^N \partial_jh_j\ | I^{\beta}(\nabla u)(t)|^2\ dxdt
\,.
\end{multline}

Thanks again to the condition \eqref{eq:mmm} the left-hand side of \eqref{eq:identity2} becomes
\begin{equation*}
\int_{\partial\Omega} \big|I^{\beta}(\partial_\nu u )(t)-I^{\beta}(\partial_\nu u )(\tau)\big|^2d\sigma 
\,.
\end{equation*}
Therefore, if we multiple both terms of \eqref{eq:identity2} by $\frac1{|t-\tau|^{1+2\beta}}$ and then integrate over $[0,T]\times[0,T]$, we have
 \begin{multline}\label{eq:identity21}
\big[ I^{\beta}(\partial_\nu u) \big]_{H^{\beta}(0,T;L^2(\partial\Omega))}^2
\\
=
2\int_0^T\int_0^T\frac1{|t-\tau|^{1+2\beta}}\langle I^{\beta}(\partial_t^{\alpha}u)(t)-I^{\beta}(\partial_{t}^{\alpha}u)(\tau),h\cdot \big(I^{\beta}(\nabla u)(t)-I^{\beta}(\nabla u)(\tau)\big)\rangle_{-\theta,\theta}\ dt d\tau
\\
+2\int_0^T\int_0^T\frac1{|t-\tau|^{1+2\beta}}
\sum_{i,j=1}^N\int_\Omega
\partial_i  h_j  \big(I^{\beta}(\partial_i u)(t)-I^{\beta}(\partial_i u)(\tau)\big)  \big(I^{\beta}(\partial_j u)(t)-I^{\beta}(\partial_j u)(\tau)\big)\ dx \ dt d\tau
\\
-
\int_0^T\int_0^T\frac1{|t-\tau|^{1+2\beta}}\int_{\Omega} \sum_{j=1}^N\partial_jh_j\ |I^{\beta}(\nabla u )(t)-I^{\beta}(\nabla u )(\tau)|^2\ dx \ dt d\tau
\,.
\end{multline}
To estimate the first term on the right-hand side of the above identity, we note that
\begin{multline}\label{eq:u-tt1}
2\int_0^T\int_0^T\frac1{|t-\tau|^{1+2\beta}}\langle I^{\beta}(\partial_t^{\alpha}u)(t)-I^{\beta}(\partial_{t}^{\alpha}u)(\tau),h\cdot 
\big(I^{\beta}(\nabla u)(t)-I^{\beta}(\nabla u)(\tau)\big)\rangle_{-\theta,\theta}\ dt d\tau
\\
\le
C\big[I^{\beta}(\partial_t^{\alpha}u)\big]_{H^{\beta}(0,T;D(A^{-\theta}))}^2
+C\big[I^{\beta}(\nabla u)\big]_{H^{\beta}(0,T;D(A^{\theta}))}^2
\,.
\end{multline}
If we choose $\theta\in \big(\frac{\alpha-1}{2\alpha},\frac1{2\alpha}\big)$, then we can apply Theorem \ref{th:reg-l2} to get $\partial_t^{\alpha}u\in L^2(0,T;D(A^{-\theta}))$ and 
$\nabla u\in L^2(0,T;D(A^{\theta}))$. Therefore, thanks to Theorem \ref{th:R-Lop} we have
\begin{equation}\label{eq:regular1}
\begin{split}
\big\|I^{\beta}(\partial_t^{\alpha}u)\big\|_{H^{\beta}(0,T;D(A^{-\theta}))}
&\sim
\|\partial_t^{\alpha}u\|_{L^2(0,T;D(A^{-\theta}))}
\,,
\\
\big\|I^{\beta}(\nabla u)\big\|_{H^{\beta}(0,T;D(A^{\theta}))}
&\sim
\|\nabla u\|_{L^2(0,T;D(A^{\theta}))}
\,.
\end{split}
\end{equation}

\end{Proof}

\end{document}